\documentclass[number,seceqn,dvips]{arxbj}
\usepackage{subenv}
\usepackage{cursive}  

\aid{0}
\volume{16}
\issue{4}
\pubyear{2010}
\firstpage{1191}
\lastpage{1207}
\doi{10.3150/09-BEJ237}

\makeatletter
\newtheorem{lem}{Lemma}[section]
\renewcommand{\emptyset}{\varnothing}
\makeatother

\begin{document}
\begin{frontmatter}

\title{Multivariate saddlepoint approximations in tail probability and
conditional inference}
\runtitle{Multivariate saddlepoint approximations}

\begin{aug}
\author{\fnms{John} \snm{Kolassa}\corref{}\thanksref{e1}\ead[label=e1,mark]{kolassa@stat.rutgers.edu}} \and
\author{\fnms{Jixin} \snm{Li}\thanksref{e2}\ead[label=e2,mark]{jixli@stat.rutgers.edu}}
\runauthor{J. Kolassa and J. Li}
\address{Department of Statistics, Rutgers University, Hill Center, Busch
Campus, 110 Frelinghuysen Road, Piscataway, NJ 08854-8019, USA. \printead{e1}; \printead*{e2}}
\end{aug}

\received{\smonth{9} \syear{2008}}
\revised{\smonth{8} \syear{2009}}

%
\begin{abstract}
We extend known saddlepoint tail probability approximations to
multivariate cases, including multivariate conditional cases. Our
approximation applies to both continuous and lattice variables, and
requires the existence of a cumulant generating function. The method
is applied to some examples, including a real data set from a
case-control study of endometrial cancer. The method contains less
terms and is easier to implement than existing methods, while showing an
accuracy comparable to those methods.
\end{abstract}

%

\end{frontmatter}

\section{Introduction}

Let $\mathbf{X}_1,\mathbf{X}_2,\ldots,\mathbf{X}_n$ be independent
and identically distributed random vectors from a density
$f_{\mathbf{X}}(\cdot)$ on $\mathbf{R}^d$. We construct an accurate
multivariate saddlepoint approximation of the tail probability of
the mean random vector
$\bar{\mathbf{X}}=(\mathbf{X}_1+\mathbf{X}_2+\cdots+\mathbf{X}_n)/n$.
We also develop a similar approximation for conditional tail
probabilities. The approximation has a relative error of $\mathrm{O}(n^{-1})$,
uniformly over a compact set of $\bar{\mathbf{x}}$, a realization of
$\bar{\mathbf{X}}$, under some general conditions. Our method
utilizes the likelihood ratio statistic, routinely calculated by
standard software, which makes the approximation easy to implement.

The Edgeworth expansion is a natural competitor to the saddlepoint
approximation. This expansion has a uniformly bounded absolute error
and works well in the center of the distribution being approximated.
However, the approximation deteriorates at the far tail of the
distribution, where it can sometimes even attain negative values.
\cite{Daniels1954} first applied saddlepoint techniques to the
approximation of a probability density function. Saddlepoint
approximation addresses the problem of degradation outside a region
of radius $\mathrm{O}(n^{-1/2})$ about $E(\mathbf{X}_i)$ by bounding
the relative error, rather than the absolute error, of the
approximation over the admissible range.

\cite{Daniels1954} discussed approximating the density of $\bar{X}$
when the dimension $d=1$, that is, the univariate case. The
approximation achieved a relative error of $\mathrm{O}(n^{-1})$ uniformly
over the whole admissible range of the variable, under some
conditions. The method uses the Fourier inversion formula, which
involves moment generating, or characteristic, functions and complex
integration. In this approach, the path of integration is shifted
so that it passes through the saddlepoint of the integrand and
follows the steepest descent curve at the neighborhood of the
saddlepoint. The asymptotic property follows from a lemma due to
\cite{Watson1948}.

Extensions of univariate saddlepoint approximation of tail
probabilities $P(\bar{X}>\bar{x})$ for the means of independent
random variables have also been studied. This calculation is more
difficult, in that, unlike the density function case, the integrand
of the Fourier inversion integral for tail probabilities has a pole
at zero.

\cite{Robinson1982} presented a general saddlepoint approximation
technique that can be applied to tail probability approximation,
based on Laplace approximation of the integrated saddlepoint
density, with an error of $\mathrm{O}(n^{-1})$. Robinson used an argument
involving a conjugate exponentially shifted distribution family and
the Edgeworth expansion. The terms of the expansion can then be
integrated termwise. There is no direct explicit formula for the
integration of each term, but the terms may be computed recursively.
This method applies when $\bar{x}\geq E(X)$. When $\bar{x}<E(X)$,
Boole's law and reflection of the distribution must be used.

\cite{LugannaniRice1980} provided an alternative approximation.
\cite{Daniels1987} derived this technique, using a transformation of
variables to directly address the local quadratic behavior of the
numerator exponent. The integral is then split into two parts, one
which contains a pole, but can be integrated exactly and explicitly,
and the other which only has removable singularities and can be
expanded and approximated accurately. The virtue of this method is
that the approximation is compact and can be computed without
recursion, and the formula is valid over the whole range of
admissible $\bar{x}$.

\cite{Reid1988} thoroughly discussed the usefulness of the saddlepoint
method in a review of the method focusing on a variety
of applications to statistical inference.

\cite{Kolassa2003} generalized the univariate Robinson approach
under the Daniels framework and achieved an error of size
$\mathrm{O}(n^{-1})$. The method uses integral expressions for the tail
probability in the multivariate case and presents a multivariate
expansion of the numerator of the integrand and a termwise
multivariate integration using recursion. This approach shares the
drawback of Robinson's approach in that it requires a positivity
constraint on the ordinate.

\cite{Wang1991} generalized Lugannani and Rice's method to the case
of a bivariate probability distribution function using variable
transformations. \cite{Kolassa2003} used a different method of
proof and showed that the error term is of order $\mathrm{O}(n^{-1})$; his
method is limited to $d=2$. Furthermore, Wang's development involves
an inversion integral in which the pole of one variable depends on
the values of other variables in a fundamentally nonlinear way.

Wang's proof of the error rate in the neighborhood of the pole is
incomplete. In this paper, a~way of effectively extending
Lugannani and Rice's method to the multivariate case, which uses a
different transformation formula from Wang's and can be used in the
case $d>2$, is proposed. The method uses fewer terms and can be extended
to multivariate conditional cases.

Our proposed saddlepoint approximation may be used to test null and
alternative hypotheses concerning a multivariate parameter when the
hypotheses are specified by systems of linear inequalities.
\cite{Kolassa2004} applied the method of \cite{Kolassa2003}, in
conjunction with the adjusted profile likelihood, in such a case.
For instance, \cite{Kolassa2004} refers to data presented by
\cite{Stokes1995} on 63 case-control pairs of women with endometrial
cancer. The occurrence of endometrial cancer is influenced by
explanatory variables including gall bladder disease, hypertension
and non-estrogen drug use. The test of whether hypertension or
non-estrogen drug use is associated with an increase in endometrial
cancer will be performed, conditional on the sufficient statistic
value associated with gall bladder disease.

The remainder of the paper is organized as follows. Section~\ref{sec2}
provides the unified framework under which both unconditional and
conditional tail probability approximations are considered. Section~\ref{sec3} derives formulas for multivariate unconditional distributions.
Section~\ref{sec4} focuses on conditional distributions. Section~\ref{sec5} presents
five examples and shows the approximation results.

\section{Multivariate extension}\label{sec2}

The unconditional and conditional tail probability approximation
share some common characteristics. We derive them in a unified way.
Applying the Fourier inversion theorem and Fubini's theorem, as in
\cite{Kolassa2003}, we find that both the unconditional and
conditional tail probability approximations require evaluation of an
integral of the form
%
\begin{equation}
\label{root_form}
\frac{n^{d-d_0}}{(2\curpi
\mathrm{i})^d}\int_{\mathbf{c}-\mathrm{i}\mathbf{K}}^{\mathbf{c}+\mathrm{i}\mathbf{K}}\frac
{\exp(n[K(\bolds{\tau})-\bolds{\tau}^\mathrm{T}\mathbf{t}^*])}{\prod
_{j=1}^{d_0}\rho(\tau_j)}\,\mathrm{d}\bolds{\tau} ,
\end{equation}
where $K$ is the cumulant generating function, which is the natural
logarithm of the moment generating function, and $\mathbf{c}$ is any
positive $d$-dimensional vector. This will be discussed in Section~\ref{sec4}.
In the unconditional case, for continuous variables, $\mathbf{K}$
is a vector of length $d$, with every entry infinity,
$\mathbf{t}^*=\mathbf{t}$ and $\rho(\tau)=\tau$; for unit lattice,
$\mathbf{K}$ is a vector of length $d$, with every entry $\curpi$,
$\mathbf{t}^*$ is~$\mathbf{t}$ corrected for continuity,
$\rho(\tau)=2\sinh(\tau/2)$ and $d=d_0$. In the conditional case,
the setting is the same, except that $d_0$ equals $d$ minus the
dimension of the conditioning variables.

\cite{Daniels1987} recast a great deal of the saddlepoint literature
in terms of inversion integrals of the form~(\ref{root_form}), rescaled
so that the exponent is exactly quadratic. This rescaling includes
the multiplier for the linear term in the exponent; this linear term
is the signed root of the likelihood ratio statistic. The idea of
using the modified signed likelihood ratio statistic was proposed in
\cite{Jensen1992}. \cite{Kolassa1997} defines a multivariate version
of this reparameterization and also defines the multiplier for the
linear terms; again, these are signed roots of likelihood ratio
statistics, but, this time, for a sequence of nested models:
\[
-\frac{1}{2}\hat{\mathbf{w}}^\mathrm{T}\hat{\mathbf{w}}=\min_{\bolds{\gamma}}
\bigl(K(\bolds{\gamma})-\bolds{\gamma}^\mathrm{T}\mathbf{t}^*\bigr)
\]
and
\[
-\frac{1}{2}(\mathbf{w}-\hat{\mathbf{w}})^\mathrm{T}(\mathbf{w}-
\hat{\mathbf{w}})=K(\bolds{\tau})-\bolds{\tau}^\mathrm{T}\mathbf{t}^*-\min_{\bolds
{\gamma}}
\bigl(K(\bolds{\gamma})-\bolds{\gamma}^\mathrm{T}\mathbf{t}^*\bigr) .
\]

Further specification of $\hat{\mathbf{w}}$ and $\mathbf{w}$ is
needed. For any vector $\mathbf{v}$ of length $d$, let
$\mathbf{v}_j$ be the vector consisting of the first $j$ elements,
that is, $(v_1,v_2,\ldots,v_j)^\mathrm{T}$. For instance,\vspace*{-1pt}
$\bolds{\gamma}_j=(\gamma_1,\gamma_2,\ldots,\gamma_j)^\mathrm{T}$,
$\bolds{\tau}_j=(\tau_1,\tau_2,\ldots,\tau_j)^\mathrm{T}$ and $\mathbf{0}_j$ is
the zero vector $(0,0,\ldots,0)^\mathrm{T}$ with dimension $j$. Let
$\mathbf{v}_{-j}$ be the vector consisting all but the first $j$
elements of $\mathbf{v}$, that is, $(v_{j+1},v_{j+2},\ldots,v_d)^\mathrm{T}$.
\cite{Kolassa1997}, Chapter 6 defines $\hat{\mathbf{w}}$ and
$\mathbf{w}$ using
%
\begin{subequation}
\begin{eqnarray}
-\frac{1}{2}\hat{w}_j^2&=&\min_{\bolds{\gamma},
\bolds{\gamma}_{j-1}=\mathbf{0}_{j-1}}
\bigl(K(\bolds{\gamma})-\bolds{\gamma}^\mathrm{T}\mathbf{t}^*\bigr)- \min_{\bolds{\gamma},
\bolds{\gamma}_j=\mathbf{0}_j} \bigl(K(\bolds{\gamma})-\bolds{\gamma}^\mathrm{T}\mathbf
{t}^*\bigr) ,\label{def_w_hat_min}\\
-\frac{1}{2}(w_j-\hat{w}_j)^2&=&\min_{\bolds{\gamma},
\bolds{\gamma}_{j-1}=\bolds{\tau}_{j-1}}
\bigl(K(\bolds{\gamma})-\bolds{\gamma}^\mathrm{T}\mathbf{t}^*\bigr)- \min_{\bolds{\gamma},
\bolds{\gamma}_j=\bolds{\tau}_j}
\bigl(K(\bolds{\gamma})-\bolds{\gamma}^\mathrm{T}\mathbf{t}^*\bigr)\label{def_w_min} .
\end{eqnarray}
\end{subequation}

This definition is not invariant with regard to the order of the
coordinates. Also, note that $w_j$ is a function of only
$\bolds{\tau}_j$, but not of any element of $\bolds{\tau}_{-j}$ $\forall{j}$.
The same holds true for $\tau_j$ as a function of
$\mathbf{w}$.

We now construct more explicit formulas for $\hat{\mathbf{w}}$ and
$\mathbf{w}$. Let
\[
\tilde{\bolds{\tau}}_j(\bolds{\gamma}_j)=(\gamma_1,\gamma_2,\ldots
,\gamma_j,\tilde{\tau}_{j+1}(\bolds{\gamma}_j),
\tilde{\tau}_{j+2}(\bolds{\gamma}_j),\ldots,\tilde{\tau}_{d}(\bolds
{\gamma}_j))
\]
be the minimizer of $(K(\bolds{\gamma})-\bolds{\gamma}^\mathrm{T}\mathbf{t}^*)$
when the first $j$ variables are fixed. The function
$\tilde{\tau}_{k}(\bolds{\gamma}_j)$ above is the minimizer for
variable $k$ when the first $j$ variables are fixed, for $k>j$.

Using the above notation, the definitions of $\hat{\mathbf{w}}$ and
$\mathbf{w}$ can be rewritten as
%
\begin{subequation}
\begin{eqnarray}
-\frac{1}{2}\hat{w}_j^2&=&K(\tilde{\bolds{\tau}}_{j-1}(\mathbf
{0}_{j-1}))-\tilde{\bolds{\tau}}_{j-1}(\mathbf{0}_{j-1})^\mathrm{T}\mathbf{t}^*
-
\bigl(K(\tilde{\bolds{\tau}}_j(\mathbf{0}_j))-\tilde{\bolds{\tau
}}_j(\mathbf{0}_j)^\mathrm{T}\mathbf{t}^*\bigr) ,\qquad \label{def_w_hat_squared}\\
-\frac{1}{2}(w_j-\hat{w}_j)^2&=&K(\tilde{\bolds{\tau}}_{j-1}(\bolds{\tau
}_{j-1}))-\tilde{\bolds{\tau}}_{j-1}(\bolds{\tau}_{j-1})^\mathrm{T}\mathbf{t}^*
-
\bigl(K(\tilde{\bolds{\tau}}_j(\bolds{\tau}_j))-\tilde{\bolds{\tau}}_j(\bolds
{\tau}_j)^\mathrm{T}\mathbf{t}^*\bigr)\label{def_w_squared} ,\qquad
\end{eqnarray}
\end{subequation}
where $\tilde{\bolds{\tau}}_{j-1}(\cdot)$ is set to $\hat{\bolds{\tau}}$
when $j=1$ for succinctness of expression.

By choosing a sign to make $\hat{\mathbf{w}}$ and $\mathbf{w}$
increasing functions of $\hat{\bolds{\tau}}$ and $\bolds{\tau}$, we can
further specify them as follows:
\begin{subequation}
\begin{eqnarray}
\hat{w}_j&=& \operatorname{sign} (\tilde{\tau}_{j}(\mathbf{0}_{j-1}))\quad
\label{def_w_hat_signed}\nonumber\\[-8pt]\\[-8pt]
&&{}\times\sqrt{-2\bigl[K(\tilde{\bolds{\tau}}_{j-1}(\mathbf{0}_{j-1}))-\tilde{\bolds
{\tau}}_{j-1}(\mathbf{0}_{j-1})^\mathrm{T}\mathbf{t}^*
-
\bigl(K(\tilde{\bolds{\tau}}_j(\mathbf{0}_j))-\tilde{\bolds{\tau
}}_j(\mathbf{0}_j)^\mathrm{T}\mathbf{t}^*\bigr)\bigr]} ,\quad\nonumber\\
w_j&=&\hat{w}_j+\operatorname{sign}\bigl(\tau_j-\tilde{\tau}_{j}(\bolds{\tau
}_{j-1})\bigr)\quad\label{def_w_signed}\nonumber\\[-8pt]\\[-8pt]
&&{}\times\sqrt{-2\bigl[K(\tilde{\bolds{\tau}}_{j-1}(\bolds{\tau}_{j-1}))-\tilde{\bolds
{\tau}}_{j-1}(\bolds{\tau}_{j-1})^\mathrm{T}\mathbf{t}^*
-
\bigl(K(\tilde{\bolds{\tau}}_j(\bolds{\tau}_j))-\tilde{\bolds{\tau}}_j(\bolds
{\tau}_j)^\mathrm{T}\mathbf{t}^*\bigr)\bigr]} .\quad\nonumber
\end{eqnarray}
\end{subequation}

The derivation of the \cite{LugannaniRice1980} approximation
provided by \cite{Daniels1987} requires identification of the simple
pole in the inversion integrand. We need to match zeros in the
denominator of the multivariate integrand with functions of the
variables in the new parameterization; the points at which this
matching occurs will be denoted by a tilde. The quantities above,
such as $\hat{\bolds{\tau}}$, $\hat{\mathbf{w}}$,
$\tilde{\tau}_{j}(\bolds{\tau}_{j-1})$ and functional relationships
between $\bolds{\tau}$ and $\mathbf{w}$, etcetera, can be solved
numerically by Newton--Raphson methods, or even analytically in some
cases. Finally, we define a function
$\tilde{w}_j(\mathbf{w}_{j-1})$ such that
$\tau_j(w_1,w_2,\ldots,\tilde{w}_j(\mathbf{w}_{j-1}))=0$ for $j>1$.

It can be verified that the following properties hold:
%
\begin{subequation}
\begin{eqnarray}
&\displaystyle\bolds{\tau}_j=0 \quad\mbox{if and only if}\quad \mathbf{w}_j=0 ;& \label{property_a}\\
&\displaystyle\tilde{w}_j(\mathbf{0}_{j-1})=0 \qquad\mbox{for } j>1;& \label{property_b}\\
&\displaystyle\tau_j=\tilde{\tau}_{j}(\bolds{\tau}_{j-1}) \quad\mbox{if and only if}\quad
w_j=\hat{w}_j \qquad\mbox{for } j>1;& \label{property_c}\\
&\displaystyle\bolds{\tau}_j=\hat{\bolds{\tau}}_j \quad\mbox{if and only if}\quad
\mathbf{w}_j=\hat{\mathbf{w}}_j.& \label{property_d}
\end{eqnarray}
\end{subequation}
Below, the superscript of a function denotes differentiation with
respect to the corresponding argument of the function. We will
employ the same use of superscripts in the subsequent text of the
paper, except that when the superscript is a set, it denotes
difference, as defined at the end of this section. Also, a
superscripted ``T'' denotes
the transpose of matrix. We can obtain
${\check{w}}_j=\tilde{w}_j(\hat{\mathbf{w}}_{j-1})$ and
${\check{w}}_j^k=\tilde{w}_j^k(\hat{\mathbf{w}}_{j-1})$, which will
be used in later sections. Substituting $w_j=\check{w}_j$, $\tau_j=0$,
$\bolds{\tau}_{j-1}=\hat{\bolds{\tau}}_{j-1}$ and
$\bolds{\tau}_j=(\hat{\tau}_1,\hat{\tau}_2,\ldots,\hat{\tau
}_{j-1},0)^\mathrm{T}=(\hat{\bolds{\tau}}_{j-1},0)^\mathrm{T}$
into~(\ref{def_w_signed}), we obtain
%
\begin{equation}
\label{w_check}
{\check{w}}_j=\hat{w}_j+\operatorname{sign}(0-\hat{\tau}_j)
\sqrt
{-2\bigl[K(\hat{\bolds{\tau}})-\hat{\bolds{\tau}}^\mathrm{T}\mathbf{t}^*
-
\bigl(K(\tilde{\bolds{\tau}}_j(\hat{\bolds{\tau}}_{j-1},0))-(\hat{\bolds
{\tau}}_{j-1},0)^\mathrm{T}\mathbf{t}^*\bigr)\bigr]}.
\end{equation}
Differentiating (\ref{def_w_squared}) with respect to $w_k$ and
rearranging terms, we obtain
%
\begin{equation}
\label{w_check_deriv}
{\check{w}}_j^k=\sum_{l=k}^{j-1}\biggl({K^l(\tilde{\bolds{\tau
}}_j(\hat{\bolds{\tau}}_{j-1},0))
\cdot\frac{\mathrm{d}\tau_l}{\mathrm{d}w_k}\bigg |_{\hat{\mathbf{w}}_l}}-t_l^*\biggr)\big/(\check
{w}_j-\hat{w}_j)\\
\end{equation}
for $k<j$. The derivatives $\frac{\mathrm{d}\tau_l}{\mathrm{d}w_k}$ evaluated at the point
$\mathbf{w}_l$ can be obtained by differentiating
(\ref{def_w_squared}) with respect to $w_k$ once or twice, depending
on whether or not $w_j=\hat{w}_j$, and solving the resulting
system of equations. In particular, we are interested in
%
\begin{equation}
\label{deriv_j} \frac{\mathrm{d}\tau_j}{\mathrm{d}w_j} \bigg|_{\mathbf{w}_j}=
\cases{\displaystyle
\sqrt{\frac{1}{\sum_{l=j}^d
K^{jl}(\tilde{\bolds{\tau}}_{j-1}(\bolds{\tau}_{j-1}))\tau_l^j(
[\tilde{\bolds{\tau}}_{j-1}(\bolds{\tau}_{j-1})]_j)}}
, &\quad if $w_j=\hat{w}_j$, \cr\displaystyle
\frac{w_j-\hat{w}_j}{K^j(\tilde{\bolds{\tau}}_j(\tau_j))-t_j^*}
,&\quad if
$w_j \ne\hat{w}_j$,
}
\end{equation}
for $j\leq d_0$, where $[\cdot]_j$ denotes the first $j$
elements, and
%
\begin{equation}
\label{prod_deriv}
\prod_{j=d_0+1}^d\frac{\mathrm{d}\tau_j}{\mathrm{d}w_j} \Bigg|_{(\mathbf{w}_{d_0},
\hat{\mathbf{w}}_{-d_0})}=
\prod_{j=d_0+1}^d\sqrt{\frac{1}{\sum_{l=j}^d
K^{jl}(\tilde{\bolds{\tau}}_{d_0}(\bolds{\tau}_{d_0}))\tau_l^j(
[\tilde{\bolds{\tau}}_{d_0}(\bolds{\tau}_{d_0})]_j)}},
\end{equation}
where, for succinctness of expression, we define $\tau_l^j(\cdot)$ to be
$1$ when $l=j$. For $l>j$, we obtain $\tau_l^j(\cdot)$ by
differentiating both sides of the definition of $\tau_l^j(\cdot)$,
that is, $K^l(\cdot)=t_l^*$ with respect to $\tau_j$ $\forall l>j$, and
solving the system of equations.

Under this transformation of variables from $\bolds{\tau}$ to $\mathbf{w}$,
the Jacobian is just the product of the diagonal terms of the
Jacobian matrix and (\ref{root_form}) can be expressed as
%
\begin{eqnarray}
\label{form_G}
\hspace*{-36pt}&&\frac{n^{d-d_0}}{(2\curpi\mathrm{i})^d}\int_{\hat{\mathbf{w}}-\mathrm{i}\mathbf{K}}^{\hat{\mathbf
{w}}+\mathrm{i}\mathbf{K}}
\frac{\exp(n[(1/2)\mathbf{w}^\mathrm{T}\mathbf{w}-\hat{\mathbf
{w}}^\mathrm{T}\mathbf{w}])}
{\prod_{j=1}^{d_0}\rho(\tau_j(\mathbf{w}_j))}\prod_{j=1}^d{\frac
{\mathrm{d}\tau_j}{\mathrm{d}w_j}}\,\mathrm{d}\mathbf{w}\nonumber\\
\hspace*{-36pt}&&\quad=
\frac{n^{d-d_0}}{(2\curpi\mathrm{i})^d}\int_{\hat{\mathbf{w}}-\mathrm{i}\mathbf{K}}^{\hat{\mathbf
{w}}+\mathrm{i}\mathbf{K}}
\frac{\exp(n[(1/2)\mathbf{w}^\mathrm{T}\mathbf{w}-\hat{\mathbf
{w}}^\mathrm{T}\mathbf{w}])}
{\prod_{j=1}^{d_0}{(w_j-\tilde{w}_j(\mathbf{w}_{j-1}))}}\cdot
\prod_{j=1}^d
{\frac{\mathrm{d}\tau_j}{\mathrm{d}w_j}}\frac{\prod_{j=1}^{d_0}{(w_j-\tilde
{w}_j(\mathbf{w}_{j-1}))}}
{\prod_{j=1}^{d_0}\rho(\tau_j(\mathbf{w}_j))}\,\mathrm{d}\mathbf{w}\\
\hspace*{-36pt}&&\quad\sim
\frac{n^{d-d_0}}{(2\curpi\mathrm{i})^d}\int_{\hat{\mathbf{w}}-\mathrm{i}\mathbf{K}}^{\hat{\mathbf
{w}}+\mathrm{i}\mathbf{K}}
\frac{\exp(n[(1/2)\mathbf{w}^\mathrm{T}\mathbf{w}-\hat{\mathbf
{w}}^\mathrm{T}\mathbf{w}])}
{\prod_{j=1}^{d_0}{(w_j-\tilde{w}_j(\mathbf{w}_{j-1}))}}
G(\bolds{\tau})\,\mathrm{d}\mathbf{w} ,\nonumber
\end{eqnarray}
where
\[G(\bolds{\tau})=\frac{\prod_{j=1}^{d_0}{(w_j-\tilde{w}_j(\mathbf
{w}_{j-1}))}}
{\prod_{j=1}^{d_0}\rho(\tau_j(\mathbf{w}_j))}\prod
_{j=1}^{d_0}\frac{\mathrm{d}\tau_j}{\mathrm{d}w_j}\cdot\prod_{j=d_0+1}^{d}
{\frac{\mathrm{d}\tau_j}{\mathrm{d}w_j}}\Bigg |_{(\mathbf{w}_{d_0},\hat{\mathbf{w}}_{-d_0})}
\]
and, to simplify notation, we set $\tilde{w}_j(\mathbf{w}_{j-1})$
to zero for $j=1$. For later convenience, we write $G(\bolds{\tau})$ as
a function $\bolds{\tau}$ instead of $\mathbf{w}$. The relation $\sim$
in the last step indicates exact equality in the unconditional case,
where $d=d_0$, but holds with a relative error of $\mathrm{O}(n^{-1})$ in
the conditional case, which we will discuss in Section~\ref{sec4}. Hereafter, we
use $\sim$ to denote approximation with a relative error of
$\mathrm{O}(n^{-1})$ of both the left-hand side and the tail probability, and
we use $\dot{\sim}$ ($\sim$ with a dot above it) in the case where the
right-hand side is an approximation with a relative error of
$\mathrm{O}(n^{-1/2})$ of the left-hand side.

The last integral in (\ref{form_G}) will be evaluated by splitting
it into rather simple terms involving poles and more complicated
terms involving analytic functions. We can decompose~(\ref{form_G})
into $2^{d_0}$ terms. Let $U=\{1,2,\ldots,d_0\}$ be the index set
of integers from $1$ to $d_0$. For set $s\subseteq U$, define
$G^s(\bolds{\tau})=G(\bolds{\tau}^s)$, where the vector $\bolds{\tau}^s$ is
defined by
\[
\tau_j^s=
\cases{
\tau_j ,&\quad if  $j\in s $,\cr
0 ,&\quad if  $j\notin s $.
}
\]
For example, if $d_0=3$, then
$G^{\{1,2\}}(\bolds{\tau})=G(\tau_1,\tau_2,0)$. Now, for $t\subseteq U$,
define $H^t=\break\sum_{s\subseteq t}(-1)^{|t-s|}G^s(\bolds{\tau})$, where
$|\cdot|$ denotes the cardinality, that is, the number of elements of a
set. For example,
$H^{\{1,2\}}=G^{\{1,2\}}(\bolds{\tau})-G^{\{1\}}(\bolds{\tau})-G^{\{2\}
}(\bolds{\tau})+G^\emptyset(\bolds{\tau})=G(\tau_1,\tau_2,0)-G(\tau
_1,0,0)-G(0,\tau_2,0)+G(0,0,0)$,
where $\emptyset$ denotes the empty set. We conclude that
$G(\bolds{\tau})=\sum_{t\subseteq U}H^t$. This decomposition holds by
induction on $d_0$. Noting that $\forall s \subseteq U$ and $a\in
s$, $H^s(\bolds{\tau}^{\{a\}})=0$, we see that
\[
\frac{H^t(\bolds{\tau})}{\prod_{j\in{t}}{(w_j-\tilde{w}_j(\mathbf
{w}_{j-1}))}}
\]
is analytic. In other words, $|t|$ product terms in the denominator
of the integrand in (\ref{form_G}) are `absorbed' by
$H^t(\bolds{\tau})$, leaving the remaining $(d_0-|t|)$ product terms
unabsorbed. As explained in \cite{Kolassa2003}, each term that is
absorbed contributes a relative error of $\mathrm{O}(n^{-1/2})$.
Therefore, if we let $I^t$ be the integral corresponding to $H^t$,
then we obtain
%
\begin{equation}
\label{form_decompose} \frac{n^{d-d_0}}{(2\curpi\mathrm{i})^d}
\int_{\hat{\mathbf{w}}-\mathrm{i}\mathbf{K}}^{\hat{\mathbf
{w}}+\mathrm{i}\mathbf{K}}
\frac{\exp(n[(1/2)\mathbf{w}^\mathrm{T}\mathbf{w}-\hat{\mathbf
{w}}^\mathrm{T}\mathbf{w}])}
{\prod_{j=1}^{d_0}{(w_j-\tilde{w}_j(\mathbf{w}_{j-1}))}}
G(\bolds{\tau})\,\mathrm{d}\mathbf{w}\sim\sum_{|t|\leq1,t\subseteq U}I^t.
\end{equation}
In the next two sections, we compute the $I^t$, $|t|\leq
1$, $t\subseteq U$, for distribution and conditional distribution,
respectively.\footnote{More detailed derivations and formulae for
bivariate distributions can be found at
\url{http://stat.rutgers.edu/resources/technical_reports10.html}.}

\section{Multivariate distribution approximation}\label{sec3}

In the unconditional continuous case, we have $d=d_0$ and
\[G(\bolds{\tau})=\frac{\prod_{j=1}^{d_0}{(w_j-\tilde{w}_j(\mathbf
{w}_{j-1}))}}
{\prod_{j=1}^{d_0}\tau_j(\mathbf{w}_j)}\prod_{j=1}^{d_0}\frac
{\mathrm{d}\tau_j}{\mathrm{d}w_j}.
\]
Therefore,
%
\begin{eqnarray}
I^{\emptyset}&=&\frac{1}{(2\curpi\mathrm{i})^{d_0}}\int_{\hat{\mathbf{w}}-\mathrm{i}\infty}^{\hat{\mathbf
{w}}+\mathrm{i}\infty}
\frac{\exp(n[(1/2)\mathbf{w}^\mathrm{T}\mathbf{w}-\hat{\mathbf
{w}}^\mathrm{T}\mathbf{w}])}
{\prod_{j=1}^{d_0}{(w_j-\tilde{w}_j(\mathbf{w}_{j-1}))}}
G(\mathbf{0})\,\mathrm{d}\mathbf{w}\nonumber\\[-8pt]\\[-8pt]
&=&\frac{1}{(2\curpi\mathrm{i})^{d_0}}\int_{\hat{\mathbf{w}}-\mathrm{i}\infty}^{\hat{\mathbf
{w}}+\mathrm{i}\infty}
\frac{\exp(n[(1/2)\mathbf{w}^\mathrm{T}\mathbf{w}-\hat{\mathbf
{w}}^\mathrm{T}\mathbf{w}])}
{\prod_{j=1}^{d_0}{(w_j-\tilde{w}_j(\mathbf{w}_{j-1}))}}\,\mathrm{d}\mathbf{w}\nonumber
\end{eqnarray}
since $G(\mathbf{0})=0$ by properties (\ref{property_a}) and
(\ref{property_b}).

Let $u_j=w_j-\tilde{w}_j(\mathbf{w}_{j-1})$, $\hat{\mathbf{u}}$ be
such that
$\mathbf{w}(\hat{\mathbf{u}})=\hat{\mathbf{u}}$ and
$g(\mathbf{u})=\frac{1}{2}\mathbf{w}(\mathbf{u})^\mathrm{T}\mathbf
{w}(\mathbf{u})-\hat{\mathbf{w}}^\mathrm{T}\mathbf{w}(\mathbf{u})$.
By changing variables, with Jacobian equal to 1, we have
%
\begin{equation}
\label{uncond_cont_I_empty} I^{\emptyset}=\frac{1}{(2\curpi\mathrm{i})^{d_0}}
\int_{\hat{\mathbf{u}}-\mathrm{i}\infty}^{\hat{\mathbf
{u}}+\mathrm{i}\infty}
\frac{\exp(n[g(\mathbf{u})])} {\prod_{j=1}^{d_0}{u_j}}\,\mathrm{d}\mathbf
{u} .
\end{equation}

The integration in (\ref{uncond_cont_I_empty}) cannot be integrated
out exactly in general. However, using the same argument as in
\cite{Kolassa2003}, we approximate it by expanding $g(\mathbf{u})$
about $\hat{\mathbf{u}}$ up to the third degree; after termwise
integration, the resulting approximation to $I^{\emptyset}$ has
relative error $\mathrm{O}(n^{-1})$. So, $I^{\emptyset}$ can be approximated
by
%
\begin{eqnarray}
I^{\emptyset}&=&\frac{1}{(2\curpi\mathrm{i})^{d_0}}\int_{\hat{\mathbf{u}}-\mathrm{i}\infty}^{\hat{\mathbf
{u}}+\mathrm{i}\infty}
\frac{\exp(n[\hat{g}+(1/2)\hat{g}^{jk}(u_j-\hat
{u}_j)(u_k-\hat{u}_k)])} {\prod_{j=1}^{d_0}{u_j}}\nonumber\\
&&{}\times \biggl(1+\frac{n}{6}\hat{g}^{jkl}(u_j-\hat{u}_j)(u_k-\hat
{u}_k)(u_l-\hat{u}_l) \biggr)\,\mathrm{d}\mathbf{u}\nonumber\\
&=&\frac{1}{(2\curpi\mathrm{i})^{d_0}}\int_{\hat{\mathbf{u}}-\mathrm{i}\infty}^{\hat{\mathbf
{u}}+\mathrm{i}\infty}
\frac{\exp(n[\hat{g}+(1/2)\hat{g}^{jk}(u_j-\hat
{u}_j)(u_k-\hat{u}_k)])}
{\prod_{j=1}^{d_0}{u_j}}\,\mathrm{d}\mathbf{u}\\
&&{}+\frac{1}{(2\curpi\mathrm{i})^{d_0}}\int_{\hat{\mathbf{u}}-\mathrm{i}\infty}^{\hat{\mathbf
{u}}+\mathrm{i}\infty}
\frac{\exp(n[\hat{g}+(1/2)\hat{g}^{jk}(u_j-\hat
{u}_j)(u_k-\hat{u}_k)])}
{\prod_{j=1}^{d_0}{u_j}}\nonumber\\
&&{}\times\frac{n}{6}\hat{g}^{jkl}(u_j-\hat{u}_j)(u_k-\hat{u}_k)(u_l-\hat
{u}_l)\,\mathrm{d}\mathbf{u} ,\nonumber
\end{eqnarray}
where, for brevity, we write $\hat{g}^r$ for $g^r(\hat{\mathbf{u}})$.
All  derivatives of $g$ evaluated at $\hat{\mathbf{u}}$ can
be computed and, in particular, $\hat{g}^j=0$. Here, we use tensor
notation, that is, the use of superscripts and subscripts to denote
summation over all possible combinations, by which we are able to omit the
summation symbol. The computation of the second integral is
addressed in \cite{Kolassa2003}. The details involve partial
derivatives of some functions up to the second or third degree; these
are algebraically complicated and therefore omitted here. For
the first integral, rearrange the terms in the numerator in the
order of the degree of $\mathbf{u}$. The first integral is quadratic
and can be\vspace*{-1.5pt} computed~as
%
\begin{eqnarray}
\label{uncond_cont_I_empty_result}
&&\frac{1}{(2\curpi\mathrm{i})^{d_0}}\int_{\hat{\mathbf{u}}-\mathrm{i}\infty}^{\hat{\mathbf
{u}}+\mathrm{i}\infty}
\frac{\exp(n[\hat{g}+(1/2)\hat{g}^{jk}(u_j-\hat
{u}_j)(u_k-\hat{u}_k)])}
{\prod_{j=1}^{d_0}{u_j}}\,\mathrm{d}\mathbf{u}\quad\nonumber\\
&&\quad=\frac{1}{(2\curpi\mathrm{i})^{d_0}}\int_{\hat{\mathbf{u}}-\mathrm{i}\infty}^{\hat{\mathbf
{u}}+\mathrm{i}\infty}
\frac{\exp(n[(\hat{g}+(1/2)\hat{g}^{jk}\hat{u}_j\hat
{u}_k)-\hat{g}^{jk}\hat{u}_ku_j+
(1/2)\hat{g}^{jk}u_ju_k])}
{\prod_{j=1}^{d_0}{u_j}}\,\mathrm{d}\mathbf{u}\quad\\
&&\quad=C^{\emptyset}\bar{\Phi}(\bar{\mathbf{y}}{}^{\emptyset},\bolds
{\Sigma}^{\emptyset}) ,\quad\nonumber
\end{eqnarray}
where
$C^{\emptyset}=\exp(n[\hat{g}+\frac{1}{2}\hat{g}^{jk}\hat
{u}_j\hat{u}_k)])$,
$\bar{\mathbf{y}}^{\emptyset}$ is a vector whose $j$th element is
$\sqrt{n}\hat{g}^{jk}\hat{u}_k/\sqrt{\hat{g}^{jj}}$ and $\bar
{\Phi}$
is the tail probability of a standard multivariate normal
distribution with mean 0 and covariance matrix
$\bolds{\Sigma}^{\emptyset}$ with elements
$\hat{g}^{jk}/\sqrt{\hat{g}^{jj}\hat{g}^{kk}}$. The last of the
above equations can be obtained by changing variables to
$\mathbf{v}$, where $v_j=u_j/\sqrt{\hat{g}^{jj}}$.

For $I^t$, $t={r}$,\vspace*{-1.5pt} we have
%
\begin{equation}
\label{uncond_cont_I_r}
I^{\{r\}}=\frac{1}{(2\curpi\mathrm{i})^{d_0}}\int_{\hat{\mathbf{w}}-\mathrm{i}\infty}^{\hat{\mathbf
{w}}+\mathrm{i}\infty}
\frac{\exp(n[(1/2)\mathbf{w}^\mathrm{T}\mathbf{w}-\hat{\mathbf
{w}}^\mathrm{T}\mathbf{w}])}
{\prod_{j \neq
r}(w_j-\tilde{w}(\mathbf{w}_{j-1}))}\cdot\frac{G^{\{r\}}(\bolds{\tau
})-G(\mathbf{0})}{w_r-\tilde{w}(\mathbf{w}_{r-1})}\,\mathrm{d}\mathbf{w} .
\end{equation}

We perform a similar change of variable from $\mathbf{w}$ to $\mathbf
{u}$ as
in computing $I^{\emptyset}$, except that $u_r=w_r$. We then\vspace*{-1.5pt} have
%
\begin{eqnarray}
\hspace*{-20pt}I^{\{r\}}&=&\frac{1}{(2\curpi\mathrm{i})^{d_0}}\int_{\hat{\mathbf{u}}-\mathrm{i}\infty}^{\hat{\mathbf
{u}}+\mathrm{i}\infty}
\frac{\exp(n[g^{\{r\}}(\mathbf{u})])} {\prod_{j \neq
r}u_j}h^{\{r\}}(\mathbf{u})\,\mathrm{d}\mathbf{w}\nonumber\\[-8pt]\\[-8pt]
\hspace*{-20pt}&\dot{\sim}&\frac{1}{(2\curpi\mathrm{i})^{d_0}}\int_{\hat{\mathbf{u}}-\mathrm{i}\infty}^{\hat{\mathbf
{u}}+\mathrm{i}\infty}
\frac{\exp(n[c_{00}^{\{r\}}+(1/2)(\mathbf{u}-\hat{\mathbf
{u}})^\mathrm{T}\mathbf{cc}^{\{r\}}(\mathbf{u}-\hat{\mathbf{u}})-
\mathbf{u}^{\{r\}}(\mathbf{u}-\hat{\mathbf{u}}))])} {\prod_{j \neq
r}u_j} h_{\mathbf{u}}^{\{r\}}(\mathbf{u})\,\mathrm{d}\mathbf{w} ,\nonumber
\end{eqnarray}
where $g^{\{r\}}(\mathbf{u})$ is the exponent as a function of
$\mathbf{u}$ after the change of variable,
$h_{\mathbf{u}}^{\{r\}}(\mathbf{u})=\frac{G^{\{r\}}(\bolds{\tau
})-G(\mathbf{0})}{w_r-\tilde{w}(\mathbf{w}_{r-1})}$.
$c^{\{r\}}_{00}=g^{\{r\}}(\hat{\mathbf{u}})$, $\mathbf{cc}^{\{r\}}$
is the matrix with elements
$cc_{ij}^{\{r\}}=[g^{\{r\}}]^{ij}(\hat{\mathbf{u}})$ and\vspace*{-4pt}
$\mathbf{c}^{\{r\}}$ is the vector such that
$c_{i}^{\{r\}}=[g^{\{r\}}]^{i}(\hat{\mathbf{u}})$. We can perform a
further change of variables $v_j=\sqrt{n}\sqrt{c_{jj}^{\{r\}}}u_j$
so\vspace*{-1.5pt} that
%
\begin{equation}
I^{\{r\}}\,\dot\sim\,\frac{C^{\{r\}}}{\sqrt{n}\sqrt{c_{rr}^{\{r\}
}}}\int_{\hat{\mathbf{v}}-\mathrm{i}\infty}^{\hat{\mathbf{v}}+\mathrm{i}\infty}
\frac{\exp((1/2)\mathbf{v}^\mathrm{T}\bolds{\Sigma}_{\mathbf{v}}^{\{r\}
}\mathbf{v}-\bar{\mathbf{y}}_{\mathbf{v}}^{\{r\}}\mathbf{v})}
{(2\curpi\mathrm{i})^{d_0}{\prod_{j \neq
r}v_j}}h_{\mathbf{v}}^{\{r\}}(\mathbf{v})\,\mathrm{d}\mathbf{v},
\end{equation}
where
$C^{\{r\}}=\exp(n[c^{\{r\}}_{00}+\frac{1}{2}\hat{\mathbf
{u}}^\mathrm{T}\mathbf{cc}^{\{r\}}\hat{\mathbf{u}}])$,
$\bolds{\Sigma}_{\mathbf{v}}^{\{r\}}$ is the covariance matrix with
elements
$[\bolds{\Sigma}_{\mathbf{v}}^{\{r\}}]_{ij}=c_{ij}^{\{r\}}/\sqrt
{c_{ii}^{\{r\}}c_{jj}^{\{r\}}}$
and
$[\bar{\mathbf{y}}_{\mathbf{v}}]_j^{\{r\}}=\frac{\sqrt{n}[\mathbf
{cc}^{\{r\}}\hat{\mathbf{u}}]_j}{\sqrt{c_{jj}^{\{r\}}}}$.
The function $h_{\mathbf{v}}^{\{r\}}(\mathbf{v})$ is analytic,
but $\frac{h_{\mathbf{v}}^{\{r\}}(\mathbf{v})}{\prod_{j \neq
r}v_j}$ is not analytic, and we cannot use Watson's lemma directly.
We use the following technique. Let
$t_r=[\bolds{\Sigma}_{\mathbf{v}}^{\{r\}}\mathbf{v}]_r$ and
$t_j=\sqrt{1-([\Sigma_{\mathbf{v}}^{\{r\}}]_{rj})^2}v_j$ for
$j\neq r$. Perform a change of variables to obtain
%
\begin{equation}
I^{\{r\}}\,\dot{\sim}\,\frac{C^{\{r\}}}{\sqrt{n}\sqrt{c_{11}^{\{r\}
}}}\int_{\hat{\mathbf{t}}-\mathrm{i}\infty}^{\hat{\mathbf{t}}+\mathrm{i}\infty}
\frac{\exp(Q^{\{r\}}(\mathbf{t}))} {(2\curpi\mathrm{i})^{d_0}\prod_{j\neq
r}t_j}h_{\mathbf{t}}^{\{r\}}(\mathbf{t})\,\mathrm{d}\mathbf{t },
\end{equation}
where
$Q^{\{r\}}(\mathbf{t})=\frac{1}{2}\mathbf{t}^\mathrm{T}\bolds{\Sigma}_{\mathbf
{t}}^{\{r\}}\mathbf{t}-\bar{\mathbf{y}}_{\mathbf{t}}^{\{r\}}\mathbf
{t}$, here $\bolds{\Sigma}_{\mathbf{t}}^{\{r\}}$ being the matrix with
elements $[\bolds{\Sigma}_{\mathbf{t}}^{\{r\}}]_{rj}=0$ for $j\neq r$,
\[
\bigl[\bolds{\Sigma}_{\mathbf{t}}^{\{r\}}\bigr]_{jk}=\frac{[\Sigma
_{\mathbf{v}}^{\{r\}}]_{jk}-
[\Sigma_{\mathbf{v}}^{\{r\}}]_{rj}[\Sigma_{\mathbf{v}}^{\{
r\}}]_{rk}}
{\sqrt{[\Sigma_{\mathbf{v}}^{\{r\}}]_{rj}[\Sigma_{\mathbf
{v}}^{\{r\}}]_{rk}}}
\qquad\mbox{for }j,k\neq r,
\]
 $\bar{\mathbf{y}}_{\mathbf{t}}^{\{r\}}$ being the
vector with elements
\[
\bigl[\bar{\mathbf{y}}_{\mathbf{t}}^{\{r\}}\bigr]_r=\bigl[\bar{\mathbf
{y}}_{\mathbf{v}}^{\{r\}}\bigr]_r
\quad\mbox{and}\quad
\bigl[\bar{\mathbf{y}}_{\mathbf{t}}^{\{r\}}\bigr]_j=\frac{[\bar{\mathbf
{y}}_{\mathbf{v}}^{\{r\}}]_j-
[\Sigma_{\mathbf{v}}^{\{r\}}]_{rj}[\bar{\mathbf{y}}_{\mathbf
{v}}^{\{r\}}]_r}
{\sqrt{1-([\Sigma_{\mathbf{v}}^{\{r\}}]_{rj})^2}}\qquad \mbox{for }j\neq
r.
\]

For a set $s$, let $\mathbf{t}^s$ denote the vector such that
$[\mathbf{t}^s]_k=0$ if $k\notin s$ and $[\mathbf{t}^s]_k=t_k$ if
$k\in s$. We have
%
\begin{equation}
I^{\{r\}}\,\dot{\sim}\,\frac{C^{\{r\}}}{\sqrt{n}\sqrt{c_{11}^{\{r\}
}}}\int_{\hat{\mathbf{t}}-\mathrm{i}\infty}^{\hat{\mathbf{t}}+\mathrm{i}\infty}
\frac{\exp(Q^{\{r\}}(\mathbf{t}))} {(2\curpi\mathrm{i})^{d_0}\prod_{k \neq
r}{t_k}}h_{\mathbf{t}}^{\{r\}}\bigl(\mathbf{t}^{\{r\}}\bigr)\,\mathrm{d}\mathbf{t} .
\end{equation}
The argument that the above holds follows similar reasoning as in
(\ref{form_decompose}), except that we only need to consider
the main term here. Now, because $t_r$ can
be separated after the change of variable and by Watson's lemma, we have
%
\begin{eqnarray}
\label{uncond_cont_I_r_result}
I^{\{r\}}
&\dot{\sim}&
\frac{C^{\{r\}}}{\sqrt{n}\sqrt{c_{11}^{\{r\}
}}}\int_{\hat{\mathbf{t}}^{\{r\}}-\mathrm{i}\infty}^{\hat{\mathbf{t}}^{\{
r\}}+\mathrm{i}\infty}
\frac{\exp(Q^{\{r\}}(\mathbf{t})-((1/2)t_r^2-[\bar{\mathbf
{y}}_{\mathbf{t}}^{\{r\}}]_rt_r)}
{(2\curpi\mathrm{i})^{d_0-1}\prod_{k \neq
r}{t_k}}\nonumber\\
&&{}\times
\int_{t_r-\mathrm{i}\infty}^{t_r+\mathrm{i}\infty}\frac{\exp((1/2)
t_r^2-[\bar{\mathbf{y}}_{\mathbf{t}}^{\{r\}}]_rt_r)}{2\curpi\mathrm{i}}h_{\mathbf{t}}^{\{r\}}
\bigl(\mathbf{t}^{\{r\}}\bigr)\,\mathrm{d}\mathbf{t}\\
&\dot{\sim}&
\frac{C^{\{r\}}h_{\mathbf{t}}^{\{r\}}(\hat{\mathbf
{t}}_r)}{\sqrt{nc_{11}^{\{r\}}}}\phi\bigl(\bigl[\bar{\mathbf{y}}_{\mathbf
{t}}^{\{r\}}\bigr]_r\bigr)\bar{\Phi} \bigl(\bar{\mathbf{y}}^{\{r\}},\Sigma^{\{r\}
} \bigr) ,\nonumber
\end{eqnarray}
where $\bar{\mathbf{y}}^{\{r\}}$ is
$\bar{\mathbf{y}}_{\mathbf{t}}^{\{r\}}$ with the $r$th element
removed and $\bolds{\Sigma}^{\{r\}}$ is
$\bolds{\Sigma}_{\mathbf{t}}^{\{r\}}$ with the $r$th row and column
removed.

Multivariate tail probability approximations for unit lattice
variables follow along the same lines, except that
\[
G(\bolds{\tau})=\frac{\prod_{j=1}^{d_0}{(w_j-\tilde{w}_j(\mathbf{w}_{j-1}))}}
{\prod_{j=1}^{d_0}2\sinh(\tau_j/2)(\mathbf{w}_j)}\prod
_{j=1}^{d_0}\frac{\mathrm{d}\tau_j}{\mathrm{d}w_j} .
\]
Since $\lim_{x\rightarrow0}{(2\sinh(x/2)/x)}=1$, any analytic
property in the continuous case still holds in the lattice case.

\section{Multivariate conditional distribution approximation}\label{sec4}

Consider a multivariate canonical exponential family. In practice,
we are often only interested in a subset of the parameters in a
given statistical model, with the other model parameters usually
treated as nuisance parameters. The distribution of the sufficient
statistics associated with the parameters of interest, conditional on
the sufficient statistics associated with the nuisance parameters,
contains the parameters of interest and not the nuisance parameters.
We can therefore use the conditional distributions instead of the
original distribution in the study. For instance, in testing
equality of proportions for a $2\times2$ contingency table, we
condition on the row or column margins. Another example is logistic
regression, where inference on some regression parameters is often
performed conditionally on sufficient statistics associated with
nuisance parameters.

Hypotheses involving parameters of interest may be tested by
computing the tail probabilities for the conditional distribution
$P(\mathbf{T}_{d_0}\geq
\mathbf{t}_{d_0}|\mathbf{T}_{-d_0}=\mathbf{t}_{-d_0})$.
\cite{Skovgaard1987} applies double saddlepoint approximation to the
problem in the case where $d_0=1$, $d>1$ and $\mathbf{T}$ is the mean
of independent and identically distributed random vectors. Here, we
propose a method that extends the results to $d_0>1$ and $d>d_0$,
using the idea from the previous sections.

First, consider $\mathbf{T}$, the mean of independent and
identically distributed continuous random vectors. Then
\[
P(\mathbf{T}_{d_0}\geq\mathbf{t}_{d_0}|\mathbf{T}_{-d_0}=\mathbf
{t}_{-d_0})=\frac{\int_{\mathbf{t}_{d_0}}^\infty f_{\mathbf
{T}}(y_1,\ldots,y_{d_0},t_{d_0+1},\ldots,t_d)\,\mathrm{d}\mathbf{y}_{d_0}}
{f_{\mathbf{T}_{-d_0}}(\mathbf{t}_{-d_0})} ,
\]
where
$f_{\mathbf{T}}(\cdot)$ is the joint density and
$f_{\mathbf{T}_{-d_0}}(\cdot)$ is the marginal density of
$\mathbf{T}_{-d_0}$. Again, we use the Fourier inversion formula to
obtain
%
\begin{equation}
\label{cond_tail}
P(\mathbf{T}_{d_0}\geq
\mathbf{t}_{d_0}|\mathbf{T}_{-d_0}=\mathbf{t}_{-d_0})= \frac
{n^{d-d_0}}{(2\curpi\mathrm{i})^d}\int_{\mathbf{c}-\mathrm{i}\infty}^{\mathbf{c}+\mathrm{i}\infty}
\frac{\exp(n[K(\bolds{\tau})-\bolds{\tau}^\mathrm{T}\mathbf{t}])}
{\prod_{j=1}^{d_0}\tau_j}\,\mathrm{d}\bolds{\tau}\big/f_{\mathbf{T}_{-d_0}}(\mathbf
{t}_{-d_0}) ,
\end{equation}
where $K(\bolds{\tau})$ is the cumulant generating function of the
random vector $\mathbf{T}$. The numerator is just a special case of~(\ref{root_form}).

Approximation (\ref{form_G}) holds because of the following lemma
which will allow us to apply previous unconditional results by
substituting components of $\hat{\mathbf{w}}$ for components of
$\mathbf{w}$ when the components correspond to variables in the
conditioning event.

\begin{lem}
%
\begin{eqnarray}
&&\frac{n^{d-d_0}}{(2\curpi\mathrm{i})^d}\int_{\hat{\mathbf{w}}-\mathrm{i}\mathbf{K}}^{\hat{\mathbf
{w}}+\mathrm{i}\mathbf{K}}
\frac{\exp(n[(1/2)\mathbf{w}^\mathrm{T}\mathbf{w}-\hat{\mathbf
{w}}^\mathrm{T}\mathbf{w}])}
{\prod_{j=1}^{d_0}{(w_j-\tilde{w}_j(\mathbf{w}_{j-1}))}}\cdot\prod_{j=1}^d
{\frac{\mathrm{d}\tau_j}{\mathrm{d}w_j}}\frac{\prod_{j=1}^{d_0}{(w_j-\tilde
{w}_j(\mathbf{w}_{j-1}))}}
{\prod_{j=1}^{d_0}\rho(\tau_j(\mathbf{w}_j))}\,\mathrm{d}\mathbf{w}\quad\nonumber\\[-8pt]\\[-8pt]
&&\quad{}=\frac{n^{d-d_0}}{(2\curpi\mathrm{i})^d}\int_{\hat{\mathbf{w}}-\mathrm{i}\mathbf{K}}^{\hat{\mathbf
{w}}+\mathrm{i}\mathbf{K}}
\frac{\exp(n[(1/2)\mathbf{w}^\mathrm{T}\mathbf{w}-\hat{\mathbf
{w}}^\mathrm{T}\mathbf{w}])}
{\prod_{j=1}^{d_0}{(w_j-\tilde{w}_j(\mathbf{w}_{j-1}))}}
G(\bolds{\tau})\,\mathrm{d}\mathbf{w}\,\bigl(1+\mathrm{O}(n^{-1})\bigr) ,\quad\nonumber
\end{eqnarray}
where
\[G(\bolds{\tau})=\frac{\prod_{j=1}^{d_0}{(w_j-\tilde{w}_j(\mathbf
{w}_{j-1}))}}
{\prod_{j=1}^{d_0}\rho(\tau_j(\mathbf{w}_j))}\cdot\prod_{j=1}^{d}
{\frac{\mathrm{d}\tau_j}{\mathrm{d}w_j}}\Bigg |_{(\mathbf{w}_{d_0},\hat{\mathbf
{w}}_{-d_0})} .
\]
\end{lem}
\begin{pf}
By Watson's lemma, given fixed $\mathbf{w}_{d_0}$, we have
\begin{eqnarray*}
&&\int_{\hat{\mathbf{w}}_{-d_0}-\mathrm{i}\mathbf{K}}^{\hat{\mathbf
{w}}_{-d_0}+\mathrm{i}\mathbf{K}}
\exp\biggl(n\biggl[\frac{1}{2}\mathbf{w}_{-d_0}^\mathrm{T}\mathbf{w}_{-d_0}-
\hat{\mathbf{w}}_{-d_0}^\mathrm{T}\mathbf{w}_{-d_0}\biggr]\biggr)
\prod_{j=d_0+1}^d {\frac{\mathrm{d}\tau_j}{\mathrm{d}w_j}}\,\mathrm{d}\mathbf{w}_{-d_0}\\
&&\quad=\int_{\hat{\mathbf{w}}_{-d_0}-\mathrm{i}\mathbf{K}}^{\hat{\mathbf
{w}}_{-d_0}+\mathrm{i}\mathbf{K}}
\exp\biggl(n\biggl[\frac{1}{2}\mathbf{w}_{-d_0}^\mathrm{T}\mathbf{w}_{-d_0}-
\hat{\mathbf{w}}_{-d_0}^\mathrm{T}\mathbf{w}_{-d_0}\biggr]\biggr)
\prod_{j=d_0+1}^d {\frac{\mathrm{d}\tau_j}{\mathrm{d}w_j}}
\Bigg|_{(\mathbf{w}_{d_0},\hat{\mathbf{w}}_{-d_0})}\\
&&\quad\quad{}\times \biggl(1+\frac{E(\mathbf{w}_{d_0})}{n} \biggr)\,\mathrm{d}\mathbf{w}_{-d_0} ,
\end{eqnarray*}
for some analytic function $E(\mathbf{w}_{d_0})$ of $\mathrm{O}(1)$.
Therefore,
\begin{eqnarray*}
\mathit{LHS}&=&\frac{n^{d-d_0}}{(2\curpi\mathrm{i})^d}\int_{\hat{\mathbf{w}}_{d_0}-
\mathrm{i}\mathbf{K}}^{\hat{\mathbf
{w}}_{d_0}+\mathrm{i}\mathbf{K}}
\frac{\exp(n[(1/2)\mathbf{w}_{d_0}^\mathrm{T}\mathbf{w}_{d_0}-
\hat{\mathbf{w}}_{d_0}^\mathrm{T}\mathbf{w}_{d_0}])}
{\prod_{j=1}^{d_0}\rho(\tau_j(\mathbf{w}_j))}\prod_{j=1}^{d_0}
{\frac{\mathrm{d}\tau_j}{\mathrm{d}w_j}}\\
&&{}\times\int_{\hat{\mathbf{w}}_{-d_0}-\mathrm{i}\mathbf{K}}^{\hat{\mathbf
{w}}_{-d_0}+\mathrm{i}\mathbf{K}}
\exp\biggl(n\biggl[\frac{1}{2}\mathbf{w}_{-d_0}^\mathrm{T}\mathbf{w}_{-d_0}-
\hat{\mathbf{w}}_{-d_0}^\mathrm{T}\mathbf{w}_{-d_0}\biggr]\biggr)
\prod_{j=d_0+1}^d {\frac{\mathrm{d}\tau_j}{\mathrm{d}w_j}}
\Bigg|_{(\mathbf{w}_{d_0},\hat{\mathbf{w}}_{-d_0})}\\
&&{}\times \biggl(1+\frac{E(\mathbf{w}_{d_0})}{n} \biggr)\,\mathrm{d}\mathbf{w}_{-d_0}\, \mathrm{d}\mathbf
{w}_{d_0}\\
&=&A\biggl(1+\frac{1}{n}\frac{B}{A}\biggr) ,
\end{eqnarray*}
where
\[
A=\frac{n^{d-d_0}}{(2\curpi\mathrm{i})^d}\int_{\hat{\mathbf{w}}-\mathrm{i}\mathbf{K}}^{\hat{\mathbf
{w}}+\mathrm{i}\mathbf{K}}
\frac{\exp(n[(1/2)\mathbf{w}^\mathrm{T}\mathbf{w}-\hat{\mathbf
{w}}^\mathrm{T}\mathbf{w}])}
{\prod_{j=1}^{d_0}{(w_j-\tilde{w}_j(\mathbf{w}_{j-1}))}}
G(\bolds{\tau})\,\mathrm{d}\mathbf{w}
\]
and
\[
B=\frac{n^{d-d_0}}{(2\curpi\mathrm{i})^d}\int_{\hat{\mathbf{w}}-\mathrm{i}\mathbf{K}}^{\hat{\mathbf
{w}}+\mathrm{i}\mathbf{K}}
\frac{\exp(n[(1/2)\mathbf{w}^\mathrm{T}\mathbf{w}-\hat{\mathbf
{w}}^\mathrm{T}\mathbf{w}])}
{\prod_{j=1}^{d_0}{(w_j-\tilde{w}_j(\mathbf{w}_{j-1}))}}
G(\bolds{\tau})E(\mathbf{w}_{d_0})\,\mathrm{d}\mathbf{w} .
\]
If $A$ and $B$ are expanded according to \cite{Kolassa2003}, each
integral is approximated by a tilting term times a normal
multivariate tail probability, up to relative order $\mathrm{O}(1/\sqrt{n})$.
The expression for $B$ is also multiplied by the leading term of
$E$. Hence, $A/B=\mathrm{O}(1)$ and, therefore, the left-hand side equals
$A(1+\mathrm{O}(n^{-1}))$.
\end{pf}

To deal with the denominator in (\ref{cond_tail}), \cite{Li2008}
demonstrates that
%
\begin{eqnarray}
\label{cond_tail_approx_2}
&&\hspace*{-30pt} \biggl(\frac{n}{2\curpi\mathrm{i}} \biggr)^{d-d_0}\int_{\hat{\mathbf{w}}
_{-d_0}-\mathrm{i}\infty}^{\hat{\mathbf
{w}}_{-d_0}+\mathrm{i}\infty}
\exp\biggl(n\biggl[\frac{1}{2}\mathbf{w}_{-d_0}^\mathrm{T}\mathbf{w}_{-d_0}-
\hat{\mathbf{w}}_{-d_0}^\mathrm{T}\mathbf{w}_{-d_0}\biggr]\biggr) \prod_{j=d_0+1}^d
\frac{\mathrm{d}\tau_j}{\mathrm{d}w_j}\Bigg |_{(\mathbf{0}_{d_0},\hat{\mathbf
{w}}_{-d_0})}\,\mathrm{d}\mathbf{w}_{-d_0}\nonumber\\[-8pt]\\[-8pt]
&&\hspace*{-30pt}\quad=f_{\mathbf{T}_{-d_0}}(\mathbf{t}_{-d_0})\bigl(1+\mathrm{O}(n^{-1})\bigr) .\nonumber
\end{eqnarray}
This development is similar to that of \cite{Kolassa1997}, page 147.

With continuous variables, we can decompose $A$ according to
(\ref{form_decompose}) with
%
\begin{equation}
G(\bolds{\tau})=\prod_{j=1}^{d_0} \biggl(\frac{w_j-\tilde{w}_j(\mathbf
{w}_{j-1})}{\tau_j}\frac{\mathrm{d}\tau_j}{\mathrm{d}w_j} \biggr) \prod_{j=d_0+1}^d
\frac{\mathrm{d}\tau_j}{\mathrm{d}w_j} \Bigg|_{(\mathbf{w}_{d_0},\hat{\mathbf{w}}_{-d_0})}.
\end{equation}
Denote the left-hand side of (\ref{cond_tail_approx_2}) by
$J_{-d_0}$. Note that $G(\mathbf{0})= \prod_{j=d_0+1}^d
\frac{\mathrm{d}\tau_j}{\mathrm{d}w_j} |_{(\mathbf{0}_{d_0},\hat{\mathbf{w}}_{-d_0})}$.
The main term is
then
%
\begin{eqnarray}
\label{cond_cont_I_empty_result}
\hspace*{-10pt}I^{\emptyset}&=&\int_{\hat{\mathbf{w}}-\mathrm{i}\infty}^{\hat{\mathbf
{w}}+\mathrm{i}\infty}
\frac{\exp(n[(1/2)\mathbf{w}^\mathrm{T}\mathbf{w}-\hat{\mathbf
{w}}^\mathrm{T}\hat{\mathbf{w}}])}
{(2\curpi\mathrm{i})^{d_0}\prod_{j=1}^{d_0}(w_j-\tilde{w}_j(\mathbf{w}_{j-1}))}\frac
{n^{d-d_0}}{(2
\curpi\mathrm{i})^{d-d_0}} \prod_{j=d_0+1}^d
\frac{\mathrm{d}\tau_j}{\mathrm{d}w_j} \Bigg|_{(\mathbf{0}_{d_0},\hat{\mathbf
{w}}_{-d_0})}\,\mathrm{d}\mathbf{w}\nonumber\\
\hspace*{-10pt}&=&\int_{\hat{\mathbf{w}}-\mathrm{i}\infty}^{\hat{\mathbf{w}}+\mathrm{i}\infty}
\frac{\exp(n[(1/2)\mathbf{w}^\mathrm{T}\mathbf{w}-\hat{\mathbf
{w}}^\mathrm{T}\hat{\mathbf{w}}])}
{(2\curpi\mathrm{i})^{d_0}\prod_{j=1}^{d_0}(w_j-\tilde{w}_j(\mathbf
{w}_{j-1}))}\,\mathrm{d}\mathbf{w}_{d_0}\cdot J_{-d_0}\\
\hspace*{-10pt}&\sim&\int_{\hat{\mathbf{w}}-\mathrm{i}\infty}^{\hat{\mathbf{w}}+\mathrm{i}\infty}
\frac{\exp(n[(1/2)\mathbf{w}^\mathrm{T}\mathbf{w}-\hat{\mathbf
{w}}^\mathrm{T}\hat{\mathbf{w}}])}
{(2\curpi\mathrm{i})^{d_0}\prod_{j=1}^{d_0}(w_j-\tilde{w}_j(\mathbf
{w}_{j-1}))}\,\mathrm{d}\mathbf{w}_{d_0}\cdot
f_{\mathbf{T}_{-d_0}}(\mathbf{t}_{-d_0}) ,\nonumber
\end{eqnarray}
where
%
\begin{equation}
\int_{\hat{\mathbf{w}}-\mathrm{i}\infty}^{\hat{\mathbf{w}}+\mathrm{i}\infty}
\frac{\exp(n[(1/2)\mathbf{w}^\mathrm{T}\mathbf{w}-\hat{\mathbf
{w}}^\mathrm{T}\hat{\mathbf{w}}])}
{(2\curpi\mathrm{i})^{d_0}\prod_{j=1}^{d_0}(w_j-\tilde{w}_j(\mathbf
{w}_{j-1}))}\,\mathrm{d}\mathbf{w}_{d_0}
\end{equation}
can be obtained by formula (\ref{uncond_cont_I_empty_result}).

Using the same technique as in
(\ref{uncond_cont_I_r})--(\ref{uncond_cont_I_r_result}), we have
%
\begin{eqnarray}
\label{cond_cont_I_r_result}
I^{\{r\}}&\dot{\sim}&\frac{n^{d-d_0}}{(2\curpi
d)^d}\int_{\hat{\mathbf{w}}-\mathrm{i}\infty}^{\hat{\mathbf{w}}+\mathrm{i}\infty}
\frac{\exp(n[(1/2)\mathbf{w}^\mathrm{T}\mathbf{w}-\hat{\mathbf
{w}}^\mathrm{T}\hat{\mathbf{w}}])}
{\prod_{j \neq
r}(w_j-\tilde{w}(\mathbf{w}_{j-1}))}\frac{G^{\{r\}}(\bolds{\tau
})-G(\mathbf{0})}{w_r-\tilde{w}(\mathbf{w}_{r-1})} \,\mathrm{d}\mathbf{w}\nonumber\\[-8pt]\\[-8pt]
&\dot{\sim}&\frac{C^{\{r\}}h_{\mathbf{t}}^{\{r\}}(\hat{\mathbf
{t}}_r)}{\sqrt{nc_{11}^{\{r\}}} \prod_{j=d_0+1}^d
\mathrm{d}\tau_j/\mathrm{d}w_j |_{(\mathbf{0}_{d_0},\hat{\mathbf
{w}}_{-d_0})}}\phi\bigl(\bigl[\bar{\mathbf{y}}_{\mathbf{t}}^{\{r\}}\bigr]_r\bigr)\bar
{\Phi} \bigl(\bar{\mathbf{y}}^{\{r\}},\Sigma^{\{r\}} \bigr)\cdot J_{-d_0}\nonumber
\end{eqnarray}
at $\mathrm{O}(n^{-1})$. The computation involves
$ \prod_{j=d_0+1}^d\frac{\mathrm{d}\tau_j}{\mathrm{d}w_j} |_{(w_1,w_2,\hat{\mathbf
{w}}_{-2})}$,
which can be obtained using~(\ref{prod_deriv}).

\begin{table}[b]
\tablewidth=328pt
\caption{Results of saddlepoint approximation compared with other
approximations in the continuous case}\label{tab1}
\begin{tabular*}{334pt}{@{}lllllll@{}}
\hline$\bar{y}_1$&$\bar{y}_2$&P. approx.&K. approx.&N.
approx.&Exact&Relative error\\
\hline
$2.5$&$2.5$&$9.12 \times10^{-2}$&$8.98\times10^{-2}$&$9.65\times
10^{-2}$&$9.22\times10^{-2}$&$-1.08\%$\\
$2.5$&$3.5$&$1.41\times10^{-2}$&$1.41\times10^{-2}$&$6.54\times
10^{-3}$&$1.41\times10^{-2}$&\phantom{$-$}$0.00\%$\\
$2.5$&$4.0$&$3.91\times10^{-3}$&$3.99\times10^{-3}$&$6.69\times
10^{-3}$&$3.93\times10^{-3}$&$-0.51\%$\\
$3.0$&$3.0$&$2.20\times10^{-2}$&$2.14\times10^{-2}$&$1.46\times
10^{-2}$&$2.22\times10^{-2}$&$-0.90\%$\\
$3.0$&$3.5$&$8.97\times10^{-3}$&$8.73\times10^{-3}$&$3.52\times
10^{-3}$&$8.96\times10^{-3}$&\phantom{$-$}$0.11\%$\\
$3.5$&$3.5$&$4.40\times10^{-3}$&$4.25\times10^{-3}$&$1.09\times
10^{-3}$&$4.40\times10^{-3}$&\phantom{$-$}$0.00\%$\\
$3.5$&$4.0$&$1.67\times10^{-3}$&$1.61\times10^{-3}$&$1.78\times
10^{-4}$&$1.66\times10^{-3}$&\phantom{$-$}$0.60\%$\\
$4.0$&$4.0$&$7.69\times10^{-4}$&$7.34\times10^{-4}$&$3.88\times
10^{-5}$&$7.58\times10^{-4}$&\phantom{$-$}$1.45\%$\\
\hline
\end{tabular*}
\end{table}

In summary, in the conditional case,
$P(\mathbf{T}_{d_0}>\mathbf{t}_{d_0}|\mathbf{T}_{-d_0}>\mathbf
{t}_{-d_0})\sim\sum_{|s|\leq
1,s \subseteq U}I^s/f_{-d_0}(\mathbf{t}_{-d_0})$, where
$U=\{1,2,\ldots,d_0\}$.

Similarly to the unconditional case, in the case of unit lattice
variables, we have
%
\begin{equation}
G(\bolds{\tau})=\prod_{j=1}^{d_0} \biggl(\frac{w_j-\tilde{w}_j(\mathbf
{w}_{j-1})}{2\sinh(\tau_j/2)}\frac{\mathrm{d}\tau_j}{\mathrm{d}w_j} \biggr) \prod_{j=d_0+1}^d
\frac{\mathrm{d}\tau_j}{\mathrm{d}w_j} \Bigg|_{(\mathbf{w}_{d_0},\hat{\mathbf
{w}}_{-d_0})} .
\end{equation}
Other analytic properties and formulae still hold.

\section{Five examples}\label{sec5}

We present five examples here. The fourth example is based on real
data.

In the first example, we consider the bivariate random vector
$(Y_1,Y_2)$, with $Y_1=X_1+X_2$ and $Y_2=X_2+X_3$, where $X_1$,
$X_2$ and $X_3$ are independent and identically distributed random
variables following the exponential distribution, which has a
density function $f(x)=\mathrm{e}^{-x}$ for $x>0$. The results for
approximating $P(\bar{Y}_1\geq\bar{y}_1,\bar{Y}_2\geq\bar{y}_2)$
when $n=5$ are listed in Table~\ref{tab1}, where ``P. approx.'' stands
for the saddlepoint approximation proposed in this paper, ``K. approx.''
stands for the saddlepoint approximation presented in
\cite{Kolassa2003} and ``N. approx.'' stands for bivariate normal
approximation. The ``exact'' column shows the exact tail probability
values computed in \cite{Mathematica05}. The ``relative error'' column
shows the relative error of ``P. approx.'' The results for the cases
$(\bar{y}_1,\bar{y}_2)=(2.5,3.0)$ and
$(\bar{y}_1,\bar{y}_2)=(3.0,4.0)$ are the special cases where
$\hat{w}_1=0,$ which we have mentioned, but which are omitted here because of the removable singularity.
The normal
approximation deteriorates at the far tail, while both saddlepoint
approximations show much better and more stable relative errors. In
almost all cases, the new method shows smaller relative errors than
those in \cite{Kolassa2003}.

\begin{table}
\tablewidth=327.5pt
\caption{Results of saddlepoint approximation compared with other
approximations in the unit lattice case}\label{tab2}
\begin{tabular*}{334.5pt}{@{}lllllll@{}}
\hline$\bar{y}_1$&$\bar{y}_2$&P. approx.&K. approx.&N.
approx.&Exact&Relative error\\
\hline
$4.5$&$4.5$&$1.15\times10^{-1}$&$1.16\times10^{-1}$&$1.16\times
10^{-1}$&$1.15\times10^{-1}$&\phantom{$-$}$0.00\%$\\
$4.5$&$5.0$&$4.43\times10^{-2}$&$4.51\times10^{-2}$&$4.28\times
10^{-2}$&$4.44\times10^{-2}$&$-0.23\%$\\
$4.5$&$5.5$&$1.04\times10^{-2}$&$1.05\times10^{-2}$&$8.73\times
10^{-3}$&$1.04\times10^{-2}$&\phantom{$-$}$0.00\%$\\
$4.5$&$6.0$&$1.46\times10^{-3}$&$1.45\times10^{-3}$&$9.50\times
10^{-4}$&$1.46\times10^{-3}$&\phantom{$-$}$0.00\%$\\
$5.0$&$5.0$&$2.07\times10^{-2}$&$2.12\times10^{-2}$&$1.92\times
10^{-2}$&$2.08\times10^{-2}$&$-0.48\%$\\
$5.0$&$5.5$&$5.89\times10^{-3}$&$6.04\times10^{-3}$&$4.85\times
10^{-3}$&$5.91\times10^{-3}$&$-0.34\%$\\
$5.0$&$6.0$&$9.91\times10^{-4}$&$1.01\times10^{-3}$&$6.40\times
10^{-4}$&$9.94\times10^{-4}$&$-0.30\%$\\
$5.5$&$5.5$&$2.11\times10^{-3}$&$2.16\times10^{-3}$&$1.57\times
10^{-3}$&$2.11\times10^{-3}$&\phantom{$-$}$0.00\%$\\
$5.5$&$6.0$&$4.45\times10^{-4}$&$4.56\times10^{-4}$&$2.69\times
10^{-4}$&$4.47\times10^{-4}$&$-0.45\%$\\
$6.0$&$6.0$&$1.21\times10^{-4}$&$1.24\times10^{-4}$&$6.14\times
10^{-5}$&$1.21\times10^{-4}$&\phantom{$-$}$0.00\%$\\
\hline
\end{tabular*}
\end{table}

\begin{table}[b]
\tablewidth=297pt
\caption{Results of saddlepoint approximation compared with
bivariate normal approximation in the conditional continuous case}\label{tab3}
\begin{tabular*}{303pt}{@{}lllllll@{}}
\hline$\bar{y}_1$&$\bar{y}_2$&$\bar{y}_3$&P. approx.&N.
approx.&Exact&Relative error\\
\hline
$2.0$&$2.0$&$7.0$&$4.42\times10^{-1}$&$8.04\times10^{-2}$&$4.38\times
10^{-1}$&\phantom{$-1$}$0.91\%$\\
$2.5$&$2.5$&$7.0$&$6.25\times10^{-2}$&$2.04\times10^{-2}$&$6.32\times
10^{-2}$&\phantom{$1$}$-1.11\%$\\
$2.5$&$3.0$&$7.0$&$8.00\times10^{-3}$&$4.14\times10^{-5}$&$8.54\times
10^{-3}$&\phantom{$1$}$-6.32\%$\\
$3.0$&$3.0$&$7.0$&$3.02\times10^{-4}$&$1.00\times10^{-8}$&$3.46\times
10^{-4}$&$-12.7\%$\\
$2.0$&$2.0$&$6.5$&$2.93\times10^{-1}$&$1.16\times10^{-1}$&$2.91\times
10^{-1}$&\phantom{$-1$}$0.69\%$\\
$2.0$&$3.0$&$6.5$&$1.09\times10^{-2}$&$6.48\times10^{-5}$&$1.14\times
10^{-2}$&\phantom{$1$}$-4.39\%$\\
$2.5$&$2.5$&$6.5$&$1.49\times10^{-2}$&$6.96\times10^{-4}$&$1.56\times
10^{-2}$&\phantom{$1$}$-4.49\%$\\
$2.5$&$3.0$&$6.5$&$5.25\times10^{-4}$&$1.57\times10^{-7}$&$6.09\times
10^{-4}$&$-13.8\%$\\
$3.0$&$3.0$&$6.5$&$9.63\times10^{-7}$&$3.67\times
10^{-12}$&$1.10\times10^{-6}$&\phantom{$-$}$12.5\%$\\
\hline
\end{tabular*}
\end{table}

In the second example, we consider the bivariate random vector
$(Y_1,Y_2)$, with $Y_1=X_1+X_2$ and $Y_2=X_2+X_3$, where $X_1$,
$X_2$ and $X_3$ are independent and identically distributed\vspace*{1pt} random
variables following the binomial distribution, which has a mass function
${N \choose x} p^x(1-p)^{N-x}$ for $0\leq x \leq N$. The results for
approximating $P(\bar{Y}_1\geq\bar{y}_1,\bar{Y}_2\geq\bar{y}_2)$
when $N=10$, $p=0.2$ and $n=8$ are displayed in Table~\ref{tab2}. We can
again see from the table that the normal approximation (with
adjustment for continuity) deteriorates at the far tail, while the
saddlepoint approximations show much better and more stable relative
errors. In most cases, the new approximation shows better accuracy
than that of \cite{Kolassa2003}.

The third example involves conditional distribution functions. Let
$X_i$, $i=1,2,3,$ be independent and identically distributed random
variables following the exponential distribution, as in the first
example. Consider the random vector $(Y_1,Y_2,Y_3)$ with $Y_1=X_2$,
$Y_2=X_3$ and $Y_3=X_1+X_2+X_3$. The results for approximating
$P(\bar{Y}_1\geq\bar{y}_1,\bar{Y}_2\geq\bar{y}_2|\bar{Y}_3=\bar{y}_3)$
when $n=10$ are shown below in Table~\ref{tab3}. The case where
$\bar{y}_1=2.0$, $\bar{y}_2=2.5$ and $\bar{y}_3=7.0$ is the special
case where both $\tilde{\tau}_2(0)=0$ and $\hat{w}_2=0$, as
discussed in Section~\ref{sec3}, and is omitted here. The cases where
$\bar{y}_1=2.0$, $\bar{y}_2=3.0$ and $\bar{y}_3=7.0$, and
$\bar{y}_1=2.0$, $\bar{y}_2=2.5$ and $\bar{y}_3=6.5$, are the cases
where $\hat{w}_1=0$; these are also omitted. The exact values are
computed in \cite{Mathematica05}.

\begin{table}
\tablewidth=278pt
\caption{Differences between cases and controls for endometrial
cancer data}\label{tab4}
\begin{tabular*}{278pt}{@{}lllllllll@{}}
\hline
Gall bladder disease&$-1$&$-1$&$-1$&\phantom{$-$}0&\phantom
{$-$}0&\phantom{1}0&\phantom{1}0&\phantom{1}0\\
Hypertension&$-1$&\phantom{$-$}0&\phantom{$-$}1&$-1$&$-1$&\phantom
{1}0&\phantom{1}0&\phantom{1}1\\
Non-estrogen drug use&\phantom{$-$}0&$-1$&\phantom
{$-$}0&$-1$&\phantom{$-$}0&\phantom{1}0&\phantom{1}1&\phantom{1}0\\
Number of pairs&\phantom{$-$}1&\phantom{$-$}1&\phantom
{$-$}1&\phantom{$-$}2&\phantom{$-$}6&14&10&12\\
[6pt]
Gall bladder disease&\phantom{$-$}0&\phantom{$-$}1&\phantom
{$-$}1&\phantom{$-$}1&\phantom{$-$}1&\phantom{1}1&\phantom
{1}1&\phantom{1}1\\
Hypertension&\phantom{$-$}1&$-1$&$-1$&\phantom{$-$}0&\phantom
{$-$}0&\phantom{1}0&\phantom{1}1&\phantom{1}1\\
Non-estrogen drug use&\phantom{$-$}1&\phantom{$-$}0&\phantom
{$-$}1&$-1$&\phantom{$-$}0&\phantom{1}1&\phantom{1}0&\phantom{1}1\\
Number of pairs&\phantom{$-$}4&\phantom{$-$}3&\phantom
{$-$}1&\phantom{$-$}1&\phantom{$-$}4&\phantom{1}1&\phantom
{1}1&\phantom{1}1\\
\hline
\end{tabular*}
\end{table}

\begin{table}[b]
\tablewidth=316pt
\caption{Endometrial cancer results for some $(t_2,t_3)$
instances}\label{tab5}
\begin{tabular*}{320pt}{@{}llllll@{}}
\hline
Method&$(10,13)$&$(9,12)$&$(8,11)$&$(7,10)$&$(6,9)$\\
\hline
N. app.&$3.50\times10^{-4}$&$1.78 \times10^{-3}$&$7.26\times
10^{-3}$&$2.39\times10^{-2}$&$6.39\times10^{-2}$\\
E. app.&$3.31\times10^{-4}$&$1.72 \times10^{-3}$&$7.13\times
10^{-3}$&$2.37\times10^{-2}$&$6.37\times10^{-2}$\\
K. app.&$1.51\times10^{-4}$&$1.07 \times10^{-3}$&$5.37\times
10^{-3}$&$2.01\times10^{-2}$&$5.84\times10^{-2}$\\
P. app.&$1.62\times10^{-4}$&$1.13 \times10^{-3}$&$5.60\times
10^{-3}$&$2.08\times10^{-2}$&$6.00\times10^{-2}$\\
Exact&$1.52\times10^{-4}$&$1.09 \times10^{-3}$&$5.48 \times
10^{-3}$&$2.05 \times10^{-2}$&$5.95 \times10^{-2}$\\
\hline
\end{tabular*}
\end{table}

The fourth example was used in \cite{Kolassa2003} and
\cite{Kolassa2004}, which refers to data presented in
\cite{Stokes1995}. The data consist of 63 case-control pairs of
women with endometrial cancer. The relationship between the
occurrence of endometrial cancer and explanatory variables including
gall bladder disease, hypertension and non-estrogen drug use is
modeled with logistic regression. \cite{Stokes1995} noted that the
likelihood for these data is equivalent to that of a logistic
regression in which the units of observation are the matched pairs,
the explanatory variables are those of the case member minus those
of the control member and the response variable is 1.

The number of pairs with each configuration of differences of the
three variables are shown in Table~\ref{tab4}.
Let $\mathbf{z}_j$, $j=1,2,\ldots,63$ denote the differences of
covariates between cases and controls, as given in Table~\ref{tab4}. Consider
the situation under the null hypothesis, where the linear coefficients
are zero. Let $\mathbf{Z}_j$, $j=1,2,\ldots,63,$ be the random vectors
that take value $\mathbf{z}_j$ with a probability of $\frac{1}{2}$
and $\mathbf{0}$ with a probability of $\frac{1}{2}$. Let
$\mathbf{Z}$ be matrix whose rows are $\mathbf{Z}_j$ and where
$\mathbf{T}=\mathbf{Z}'\mathbf{1}$,\vspace*{-2pt} where $\mathbf{1}$ is a column vector
with dimension 63. We then have
$K(\bolds{\tau})=\sum_jm_j[\log(\frac{1+\exp(\mathbf{z}_j\bolds{\tau
})}{2})]$.
\cite{Kolassa2004} tested the association of hypertension or
non-estrogen drug use with an increase in endometrial cancer,
conditional on the sufficient statistic value associated with gall
bladder disease. The test required evaluating the quantity
$P(T_2\geq10 \mbox{ or } T_3\geq13|T_1=9)$ for
$\mathbf{T}=(T_1,T_2,T_3)$. By Boole's law, this probability can be
computed using
\[
P(T_2 \geq10|T_1 = 9)+P(T_3 \geq13|T_1 = 9)-P(T_2 \geq10,T_3 \geq
13|T_1 = 9) .
\]
The results for approximating $P(T_2 \geq10,T_3 \geq13|T_1 = 9)$
compared to those listed in \cite{Kolassa2004} are shown in Table~\ref{tab5},
where ``N. app.'' stands for normal approximation, ``E. app.''
stands for Edgeworth approximation, ``K. app.'' stands for the
approximation presented in \cite{Kolassa2004} and ``P. app.'' is the
proposed approximation. Approximation results of $P(T_2 \geq t_2,T_3
\geq t_3|T_1 = 9)$ for other values of~$t_2$ and $t_3$ are also
listed in the table. We can see that the proposed method achieves
better results than other methods, except for the method of \cite
{Kolassa2004}, which is far more complicated
computationally.\looseness=1

\begin{table}
\tablewidth=350.5pt
\caption{Results of saddlepoint approximation compared with normal
approximations for a multivariate gamma distribution}\label{tab6}
\begin{tabular*}{356.5pt}{@{}llllllll@{}}
\hline$y_1$&$y_2$&$y_3$&P. approx.&N.
approx.&Simulation&Std. err.&Relative error\\
\hline
$5.5$&$5.5$&$5.5$&$4.93 \times10^{-2}$&$9.64\times
10^{-2}$&$5.58\times10^{-2}$&$1.42\times10^{-3}$&$-11.6\%$\\
$5.5$&$5.5$&$6.5$&$3.65 \times10^{-2}$&$7.17\times
10^{-2}$&$4.16\times10^{-2}$&$1.24\times10^{-3}$&$-12.3\%$\\
$5.5$&$6.5$&$6.5$&$2.70 \times10^{-2}$&$5.34\times
10^{-2}$&$3.26\times10^{-2}$&$1.10\times10^{-4}$&$-17.2\%$\\
$6.5$&$6.5$&$6.5$&$1.98 \times10^{-2}$&$3.99\times
10^{-2}$&$2.45\times10^{-2}$&$9.58\times10^{-4}$&$-19.2\%$\\
$6.5$&$6.5$&$7.5$&$1.45 \times10^{-2}$&$2.78\times
10^{-2}$&$1.83\times10^{-2}$&$8.31\times10^{-4}$&$-20.8\%$\\
$6.5$&$7.5$&$7.5$&$1.05 \times10^{-2}$&$1.94\times
10^{-2}$&$1.34\times10^{-2}$&$7.11\times10^{-4}$&$-21.6\%$\\
$7.5$&$7.5$&$7.5$&$0.76 \times10^{-2}$&$1.36\times
10^{-3}$&$1.04\times10^{-3}$&$1.89\times10^{-4}$&$-26.9\%$\\
\hline
\end{tabular*}
\end{table}

In the fifth example, we consider a multivariate gamma distribution,
which is the diagonal of a Wishart distribution, formed from a
3-variate normal distribution, with covariance matrix
\[
V= \pmatrix{
1 & 0.25 & 0.25 \cr
0.25 & 1 & 0.25 \cr
0.25 & 0.25 & 1}
\]
and
$n=5$. The results are listed in Table~\ref{tab6}, where ``P.
approx.'' stands for the saddlepoint approximation proposed in this
paper and ``N. approx.'' stands for bivariate normal approximation.
The ``simulation'' and ``std. err.'' column shows the simulation results
and $5\%$ standard error. The ``relative error'' column shows the
relative error of ``P. approx.'' compared with the simulation
results. We can see that the proposed approximation performs better
than the normal approximation.\looseness=1

\section*{Acknowledgements}
This research was supported in part by NSF Grant DMS-0505499. The
authors would like to thank the referees and the Associate Editor for their
helpful suggestions.

\printhistory


\begin{thebibliography}{99}

\bibitem{Daniels1954}
Daniels, H.E. (1954).
Saddlepoint approximations in statistics.
\textit{Ann. Math. Statist.} \textbf{25} 631--645.
\MR{0066602}

\bibitem{Daniels1987}
Daniels, H.E. (1987).
Tail probability approximations.
\textit{Internat. Statist. Rev.} \textbf{55} 37--46.
\MR{0962940}

\bibitem{Jensen1992}
Jensen, J.L. (1992).
The modified signed likelihood statistic and saddlepoint
approximations.
\textit{Biometrika} \textbf{79} 693--703.
\MR{1209471}

\bibitem{Kolassa1997}
Kolassa, J.E. (2006).
\textit{Series Approximation Methods in Statistics}, 3rd ed.
\textit{Lecture Notes in Statistics} \textbf{88}.
New York: Springer.
\MR{1487639}

\bibitem{Kolassa2003}
Kolassa, J.E. (2003).
Multivariate saddlepoint tail probability approximations.
\textit{Ann. Statist.} \textbf{31} 274--286.
\MR{1962507}

\bibitem{Kolassa2004}
Kolassa, J.E. (2004).
Approximate multivariate conditional inference using the adjusted
profile likelihood.
\textit{Canad. J. Statist.} \textbf{32} 5--14.
\MR{2060541}

\bibitem{Li2008}
Li, J. (2008).
Multivariate saddlepoint tail probability approximations, for
conditional and unconditional distributions, based on the signed root
of the
log likelihood ratio statistic.
Ph.D. thesis, Rutgers Univ., Dept. Statistics and
Biostatistics.

\bibitem{LugannaniRice1980}
Lugannani, R. and Rice, S. (1980).
Saddlepoint approximation for the distribution of the sum of
indepednent random variables.
\textit{Adv. in Appl. Probab.} \textbf{12} 475--490.
\MR{0569438}

\bibitem{Reid1988}
Reid, N. (1988).
Saddlepoint methods and statistical inference (with discussion).
\textit{Statist. Sci.} \textbf{3} 213--238.
\MR{0968390}

\bibitem{Robinson1982}
Robinson, J. (1982).
Saddlepoint approximation for permutation tests and confidence
intervals.
\textit{J. R. Statist. Soc. Ser. B} \textbf{44} 91--101.
\MR{0655378}

\bibitem{Skovgaard1987}
Skovgaard, I.M. (1987).
Saddlepoint expansions for conditional distributions.
\textit{J. Appl. Probab.} \textbf{24} 875--887.
\MR{0913828}

\bibitem{Stokes1995}
Stokes, M.E., David, C.S. and Koch, G.G. (1995).
\textit{Categorical Data Analysis Using the SAS System}.
Cary, NC.

\bibitem{Wang1991}
Wang, S. (1991).
Saddlepoint approximation for bivariate distribution.
\textit{J. Appl. Probab.} \textbf{27} 586--597.
\MR{1067024}

\bibitem{Watson1948}
Watson, G.N. (1948).
\textit{Theory of Bessel Functions}.
Cambridge:
Cambridge Univ. Press.

\bibitem{Mathematica05}
Wolfram Research, Inc. (2005).
\textit{Mathematica. Version 5.0}.
Champaign, IL: Wolfram Research, Inc.

\end{thebibliography}
\end{document}